# ON FAST COMPUTATION OF GRADIENTS FOR CANDECOMP/PARAFAC ALGORITHMS

ANH HUY PHAN *, PETR TICHAVSKÝ †, AND ANDRZEJ CICHOCKI ‡

**Abstract.** Product between mode-$n$ unfolding $\mathbf{Y}_{(n)}$ of an $N$-D tensor $\boldsymbol{\mathcal{Y}}$ and Khatri-Rao products of $(N-1)$ factor matrices $\mathbf{A}^{(m)}$, $m = 1, \ldots, n-1, n+1, \ldots, N$ exists in algorithms for CANDECOMP/PARAFAC (CP). If $\boldsymbol{\mathcal{Y}}$ is an error tensor of a tensor approximation, this product is the gradient of a cost function with respect to factors, and has the largest workload in most CP algorithms. In this paper, a fast method to compute this product is proposed. Experimental verification shows that the fast CP gradient can accelerate the CP_ALS algorithm 2 times and 8 times faster for factorizations of 3-D and 4-D tensors, and the speed-up ratios can be 20-30 times for higher dimensional tensors.

**1. Introduction.** Canonical polyadic decomposition also coined CANDECOMP/PARAFAC (CP) [8, 16] is a common tensor factorization which has found applications such as in chemometrics, telecommunication, analysis of fMRI data, time-varying EEG spectrum, data mining [7, 13], classification, clustering [29], stochastic PDEs [15]. For example, CP was applied to analyze the auditory tones by Carroll and Chang [8], or to vowel-sound data by Harshman [16], or to model fluorescence excitation-emission data by hidden loading components in chemometrics [3]. Applications of CP to sensor array processing and CDMA systems in telecommunications have been considered in [12, 30]. In neuroscience, Field and Graupe [14] extracted topographic components model from event-related potentials data, Mørup et al. [23] analyzed EEG data in the time-frequency domain.

Since the alternating least squares (ALS) algorithm was proposed [8, 16], there have been intensive research efforts to improve performance and accelerate convergence rate of CP algorithms. A number of particular techniques are developed such as line search extrapolation methods [2, 16, 27, 33], compression [19], or simply adding a small diagonal matrix [10]. Instead of alternating estimation, all-at-once algorithms such as the OPT algorithm [1], the PMF3, damped Gauss-Newton (dGN) algorithms [24, 33] and fast dGN [25, 26, 32] are studied to deal with problems of a slow convergence of the ALS in some cases. Another approach is to consider the CP decomposition as a joint diagonalization problem [11, 21, 22, 28, 28].

CP algorithms can speed-up convergence rate, or cope with difficult problems. However, in all CP algorithms, the largest workload is product of tensor unfoldings and all-but-one factors which has not been inadequately considered. If a tensor of size $I_1 \times I_2 \times \cdots \times I_N$ is an error tensor of a data tensor and its approximation, the products express the gradients of a cost function with respect to factors of size $I_n \times R$. Hereinafter, we call this product "CP gradient". The CP gradients with respect to all the factors have a high computational cost of order $O\left(NR \prod_{n=1}^{N} I_n\right)$. In addition, mode-$n$ tensor unfoldings with $n = 2, 3, \ldots, N-1$ are also time consuming due to accessing non-contiguous blocks of data entries and shuffling their orders stored in memory. For high dimensional data tensors such as $N \geq 4$, the CP gradients may become very computational demanding. Experimental results show that it might take several hours to factorize 7-D tensors if they comprise hundreds of millions or billions of entries (e.g., a tensor of size $10 \times 10 \times 10 \times 10 \times 10 \times 10 \times 10$) and has high rank $R = 10$.

*Brain Science Institute, RIKEN, Wakoshi, Japan (phan@brain.riken.jp).
†Institute of Information Theory and Automation, Prague, Czech (tichavsk@utia.cas.cz). The work of P. Tichavský was supported by Grant Agency of the Czech Republic through the project 102/09/1278.
‡Brain Science Institute, RIKEN, Wakoshi, Japan and Systems Research Institute, PAS, Poland (cia@brain.riken.jp).

In this paper, a fast computation method is proposed for the CP gradients. The method avoids mode-$n$ tensor unfoldings, and reduces the computational cost to $O\left(R \prod_{n=1}^{N} I_n\right)$ for computing CP gradients over all modes.

The paper is organized as follows. Notation and basic multilinear algebra are briefly reviewed in Section 2. CP model and CP gradients are shortly reviewed in this section. The fast computation method is presented in Section 3. The fast implementation of the ALS algorithm utilizing the fast CP gradient is introduced in Section 4. In Section 5 we provide examples illustrating the validity and performance of the proposed algorithm. Section 6 concludes the paper.

**2. Notation and CANDECOMP/PARAFAC (CP) model.** We shall denote a tensor by bold calligraphic letters, e.g., $\mathcal{A} \in \mathbb{R}^{I_1 \times I_2 \times \cdots \times I_N}$, matrices by bold capital letters, e.g. $\mathbf{A} = [\mathbf{a}_1, \mathbf{a}_2, \ldots, \mathbf{a}_R] \in \mathbb{R}^{I \times R}$, and vectors by bold italic letters, e.g. $\mathbf{a}_j$ or $\mathbf{I} = [I_1, I_2, \ldots, I_N]$.

An $\mathbf{i} = [i_1, i_2, \ldots, i_N]$-th entry $y_{\mathbf{i}} = \mathcal{Y}(i_1, i_2, \ldots, i_N)$ with $1 \le i_n \le I_n$, $n = 1, 2, \ldots, N$ is alternatively denoted by $y_i$ with the index $i = \text{ivec}(\mathbf{i}, \mathbf{I})$[1] defined as

$$i = \text{ivec}(\mathbf{i}, \mathbf{I}) = i_1 + \sum_{n=2}^{N} (i_n - 1) \prod_{j=1}^{n-1} I_j. \quad (2.1)$$

A vector of integer numbers is denoted by colon notation such as $\mathbf{k} = i{:}j = [i, i+1, \ldots, j-1, j]$. For example, we denote $1{:}n = [1, 2, \ldots, n]$.

Generally, we adopt notation used in [9, 20]. The Kronecker product, the Khatri-Rao (column-wise Kronecker) product, and the (element-wise) Hadamard product and division are denoted respectively by $\otimes, \odot, \circledast, \oslash$ [9, 20].

NOTATION 2.1. **(Hadamard and Khatri-Rao products of matrices)** *Given a set of N matrices $\mathbf{A}^{(n)} \in \mathbb{R}^{I_n \times R}$, $n = 1, 2, \ldots, N$, Hadamard and Khatri-Rao products among them are denoted by*

$$\underset{k \ne n}{\circledast} \mathbf{A}^{(k)} = \mathbf{A}^{(N)} \circledast \cdots \circledast \mathbf{A}^{(n+1)} \circledast \mathbf{A}^{(n-1)} \circledast \cdots \circledast \mathbf{A}^{(1)}, \qquad I_n = I, \forall n,$$

$$\underset{n=1}{\overset{N}{\odot}} \mathbf{A}^{(n)} = \mathbf{A}^{(N)} \odot \cdots \odot \mathbf{A}^{(n)} \odot \cdots \odot \mathbf{A}^{(1)}, \qquad \forall n,$$

$$\underset{k \ne n}{\odot} \mathbf{A}^{(k)} = \mathbf{A}^{(N)} \odot \cdots \odot \mathbf{A}^{(n+1)} \odot \mathbf{A}^{(n-1)} \cdots \odot \mathbf{A}^{(1)}.$$

DEFINITION 2.1 (**Reshaping**). *The reshape operator for a tensor $\mathcal{Y} \in \mathbb{R}^{I_1 \times I_2 \times \cdots \times I_N}$ to a size specified by a vector $\mathbf{L} = [L_1, L_2, \ldots, L_M]$ with $\prod_{m=1}^{M} L_m = \prod_{n=1}^{N} I_n$ returns an M-D tensor $\mathcal{X}$, such that $\text{vec}(\mathcal{Y}) = \text{vec}(\mathcal{X})$, and is expressed as*

$$\mathcal{X} = \texttt{reshape}(\mathcal{Y}, \mathbf{L}) \quad \in \mathbb{R}^{L_1 \times L_2 \times \cdots \times L_M}. \quad (2.2)$$

Reshape does not change the order of entries in its vectorization.

DEFINITION 2.2 (**Tensor unfolding [4]**). *Unfolding a tensor $\mathcal{Y} \in \mathbb{R}^{I_1 \times I_2 \times \cdots \times I_N}$ along modes $\mathbf{r} = [r_1, r_2, \ldots, r_M]$ and $\mathbf{c} = [c_1, c_2, \ldots, c_{N-M}]$ where $[\mathbf{r}, \mathbf{c}]$ is a permutation of $[1, 2, \ldots, N]$*

---

[1] ivec is the "sub2ind" Matlab function.

*aims to rearrange this tensor to be a matrix* $\mathbf{Y}_{r \times c}$ *of size* $\prod_{k=1}^{M} I_{r_k} \times \prod_{l=1}^{N-M} I_{c_l}$ *whose entries* $(j_1, j_2)$ *are given by* $\mathbf{Y}_{r \times c}(j_1, j_2) = \mathcal{Y}(i_r, i_c)$, *where* $i_r = [i_{r_1} \ldots i_{r_M}]$, $i_c = [i_{c_1} \ldots i_{c_{N-M}}]$, $j_1 = \text{ivec}(i_r, I_r)$, $j_2 = \text{ivec}(i_c, I_c)$.

REMARK 2.1.
1. *If* $c = [c_1 < c_2 < \cdots < c_{N-M}]$, *then* $\mathbf{Y}_{r \times c}$ *is simplified to* $\mathbf{Y}_{(r)}$.
2. *If* $r = n$ *and* $c = [1, \ldots, n-1, n+1, \ldots, N]$, *we have mode-n matricization* $\mathbf{Y}_{r \times c} = \mathbf{Y}_{(n)}$.
3. $\mathbf{Y}_{r \times c} = \mathbf{Y}_{c \times r}^T$.
4. *For* $r = [1, 2, \ldots, n]$, $c = [n+1, n+2, \ldots, N]$, $\forall n$, $\mathbf{Y}_{r \times c} = \mathbf{Y}_{(r)}$ *can be expressed and efficiently performed by* reshape, *that is*

$$\mathbf{Y}_{(r)} = \texttt{reshape}(\mathcal{Y}, [J_n, K_n]), \quad J_n = \prod_{k=1}^{n} I_k, \quad K_n = \prod_{k=n+1}^{N} I_k. \quad (2.3)$$

DEFINITION 2.3. **(mode-n tensor-vector product)** *The mode-n multiplication of a tensor* $\mathcal{Y} \in \mathbb{R}^{I_1 \times I_2 \times \cdots \times I_N}$ *by a vector* $a \in \mathbb{R}^{I_n}$ *returns an* $(N-1)$*-D tensor* $\mathcal{Z}$ *defined as*

$$\text{vec}(\mathcal{Z}) = \mathbf{Y}_{(n)}^T a. \quad (2.4)$$

*Symbolically, the product is denoted by*

$$\mathcal{Z} = \mathcal{Y} \bar{\times}_n a \in \mathbb{R}^{I_1 \times \cdots \times I_{n-1} \times I_{n+1} \times \cdots \times I_N}. \quad (2.5)$$

*Tensor-vector product of a tensor* $\mathcal{Y}$ *with a set of N column vectors* $\{a\} = \{a^{(1)}, a^{(2)}, \ldots, a^{(N)}\}$ *is denoted by*

$$\mathcal{Y} \bar{\times} \{a\} = \mathcal{Y} \bar{\times}_1 a^{(1)} \bar{\times}_2 a^{(2)} \cdots \bar{\times}_N a^{(N)}. \quad (2.6)$$

DEFINITION 2.4. **(CANDECOMP/PARAFAC (CP))** *Factorize a given N-th order data tensor* $\mathcal{Y} \in \mathbb{R}^{I_1 \times I_2 \times \cdots \times I_N}$ *into a set of N component matrices (factors):* $\mathbf{A}^{(n)} = [a_1^{(n)}, a_2^{(n)}, \ldots, a_R^{(n)}] \in \mathbb{R}^{I_n \times R}$, $(n = 1, 2, \ldots, N)$ *representing the common (loading) factors [8, 16, 17], that is,*

$$\mathcal{Y} \approx \sum_{r=1}^{R} a_r^{(1)} \circ a_r^{(2)} \circ \ldots \circ a_r^{(N)} = \hat{\mathcal{Y}}, \quad (2.7)$$

*where symbol "*$\circ$*" denotes outer product. Tensor* $\hat{\mathcal{Y}}$ *is an approximation of the data tensor* $\mathcal{Y}$. *Mode-n matricization of* $\mathcal{Y}$ *can be represented as:*

$$\mathbf{Y}_{(n)} \approx \mathbf{A}^{(n)} \left( \bigodot_{k \neq n} \mathbf{A}^{(k)} \right)^T.$$

**2.1. Complexity of Tensor Unfoldings.** Tensor unfoldings are to rearrange entries of tensors to be matrices. We note that entries of the tensor $\mathcal{Y}$ are stored as a long vector $\text{vec}(\mathcal{Y})$ of the size $\prod_{n=1}^{N} I_n$ in memory. From this view point, tensor unfolding is to change the order to entries in its vectorization. The more the changes of entries take place, the slower the unfolding are. Moreover, reading data (entries) stored in non-contiguous blocks will be at a slower rate than accessing data stored in a contiguous block.

The mode-1 unfolding $\mathbf{Y}_{(1)}$ comprises $J_{-1} = I_2 I_3 \cdots J_N$ column vectors which consist of $I_1$ contiguous entries of $\mathcal{Y}$

$$\mathbf{Y}_{(1)} = \begin{bmatrix} y_1 & y_{I_1+1} & \cdots & y_{(J_{-1}-1)I_1+1} \\ y_2 & y_{I_1+1} & \cdots & y_{(J_{-1}-1)I_1+1} \\ \vdots & \vdots & \ddots & \vdots \\ y_{I_1} & y_{2I_1} & \cdots & y_{J_N} \end{bmatrix} \qquad (2.8)$$

$$= \begin{bmatrix} \mathcal{Y}(1:I_1) & \mathcal{Y}(I_1+1:2I_1) & \cdots & \mathcal{Y}((J_{-1}-1)I_1+1:J_N) \end{bmatrix}. \qquad (2.9)$$

By taking into account that $\mathbf{Y}_{(N)} = \mathbf{Y}_{(1:N-1)}^T$, we in practice compute $\mathbf{Y}_{(1:N-1)}$ instead of $\mathbf{Y}_{(N)}$. $\mathbf{Y}_{(1:N-1)}$ consists of $I_N$ vectors each of which comprises $J_{N-1}$ contiguous entries given by

$$\mathbf{Y}_{(N)}^T = \mathbf{Y}_{(1:N-1)} = \begin{bmatrix} \mathcal{Y}(1:J_{N-1}) & \mathcal{Y}(J_{N-1}+1:2J_{N-1}) & \cdots & \mathcal{Y}((I_N-1)J_{N-1}+1:J_N) \end{bmatrix}. \quad (2.10)$$

In general, unfoldings $\mathbf{Y}_{(1:n)}$ do not change the order of entries of $\mathcal{Y}$

$$\mathbf{Y}_{(1:n)} = \begin{bmatrix} \mathcal{Y}(1:J_n) & \mathcal{Y}(J_n+1:2J_n) & \cdots & \mathcal{Y}((K_n-1)J_n+1:J_N) \end{bmatrix}, \quad n = 1, 2, \ldots, N. \quad (2.11)$$

Hence, they are relatively fast. We denote by $\mathcal{X}^{(m)}$, $\boldsymbol{m} = [i_{n+1}, i_{n+2}, \ldots, i_N]$, $n$-dimensional subtensors of $\mathcal{Y}$ whose each entry is given by $\mathcal{X}^{(m)}(i_1, i_2, \ldots, i_n) = \mathcal{Y}(i_1, i_2, \ldots, i_n, i_{n+1}, \ldots, i_N)$ The mode-$n$ unfolding $\mathbf{Y}_{(n)}$ of the size $I_n \times J_{-n}$, $J_{-n} = I_1 \cdots I_{n-1} I_{n+1} \cdots I_N = J_{n-1} K_n$, can be expressed as concatenation of $K_n$ mode-$n$ unfoldings of $\mathcal{X}^{(m)}$

$$\mathbf{Y}_{(n)} = \begin{bmatrix} \mathbf{X}_{(n)}^{(1)} & \cdots & \mathbf{X}_{(n)}^{(m)} & \cdots & \mathbf{X}_{(n)}^{(M)} \end{bmatrix}, \quad M = [I_{n+1}, I_{n+2}, \ldots, I_N],$$

$$\mathbf{X}_{(n)}^{(m)} = \begin{bmatrix} y_{mJ_n+1} & \cdots & y_{mJ_n+J_{n-1}} \\ y_{mJ_n+J_{n-1}+1} & \cdots & y_{mJ_n+2J_{n-1}+1} \\ \vdots & \ddots & \vdots \\ y_{mJ_n-J_{n-1}+1} & \cdots & y_{(m+1)J_n} \end{bmatrix}, \quad m = \text{ivec}(\boldsymbol{m}, \boldsymbol{M}) - 1.$$

Therefore, most entries of $\mathbf{Y}_{(n)}$ have changed their orders. This is why the mode-$n$ unfoldings $\mathbf{Y}_{(n)}$ for $1 < n < N$ are more time consuming, and relatively slower than unfoldings $\mathbf{Y}_{(1:n)}$.

**2.2. Gradients in CP Algorithms.** We consider the cost function

$$D = \|\mathcal{Y} - \widehat{\mathcal{Y}}\|_F^2, \qquad (2.12)$$

and the gradients of this cost function with respect to the factor $\mathbf{A}^{(n)}$, $n = 1, 2, \ldots, N$ are given by [25, 33]

$$\mathbf{G}^{(n)} = \mathbf{E}_{(n)} \left( \bigodot_{k \neq n} \mathbf{A}^{(k)} \right) = \mathbf{Y}_{(n)} \left( \bigodot_{k \neq n} \mathbf{A}^{(k)} \right) - \mathbf{A}^{(n)} \left( \circledast_{k \neq n} \mathbf{A}^{(k)T} \mathbf{A}^{(k)} \right) \in \mathbb{R}^{I_n \times R}, n = 1, \ldots, N, (2.13)$$

where $\mathbf{E}_{(n)}$ denotes the mode-$n$ unfolding of the error tensor $\mathcal{E} = \mathcal{Y} - \widehat{\mathcal{Y}}$. The product $\mathbf{E}_{(n)} \left( \bigodot_{k \neq n} \mathbf{A}^{(k)} \right)$ or $\mathbf{Y}_{(n)} \left( \bigodot_{k \neq n} \mathbf{A}^{(k)} \right)$ has a computational cost of order $O(R J_N)$, and is the most expensive step in CP algorithms. Indeed, the mode-$n$ unfoldings $\mathbf{Y}_{(n)}$ for $n > 1$ are time-consuming, but are not appropriately computed. The latter product $\mathbf{Y}_{(n)} \left( \bigodot_{k \neq n} \mathbf{A}^{(k)} \right)$ is more

**Algorithm 1:** Direct Computation of $\mathbf{Y}_{(n)} \left( \bigodot_{k \neq n} \mathbf{A}^{(k)} \right)$ [5, 6]

**Input**: $\mathcal{Y}$: $(I_1 \times I_2 \times \cdots \times I_N)$, $N$ matrices $\mathbf{A}^{(n)} \in \mathbb{R}^{I_n \times R}$
**Output**: $\mathbf{G}^{(n)} = \mathbf{Y}_{(n)} \left( \bigodot_{k \neq n} \mathbf{A}^{(k)} \right) : I_n \times R$
**begin**

1    $\mathcal{Y} \leftarrow \texttt{permute}(\mathcal{Y}, [n, 1{:}n-1, n+1{:}N])$      % tensor transposition
2    $\mathbf{Y}_{(n)} \leftarrow \texttt{reshape}(\mathcal{Y}, [I_n, J_{-n}])$      % tensor unfolding $\mathbf{Y}_{(n)}$
3    $\mathbf{G}^{(n)} = \mathbf{Y}_{(n)} \left( \bigodot_{k \neq n} \mathbf{A}^{(k)} \right)$

efficient than $\mathbf{E}_{(n)} \left( \bigodot_{k \neq n} \mathbf{A}^{(k)} \right)$ in the sense of computation because it does not need to construct the error tensor $\mathcal{E}$. However, since both products involve the same mathematical expression, we also call $\mathbf{Y}_{(n)} \left( \bigodot_{k \neq n} \mathbf{A}^{(k)} \right)$ the CP gradient in which $\mathcal{Y}$ is considered as an error tensor.

The CP gradients are employed in almost all CP algorithms. For example, the alternating least squares (ALS) algorithm [2, 8, 16, 30, 31] alternatively minimizes the cost function (2.12) with an update rule given by

$$\mathbf{A}^{(n)} \leftarrow \mathbf{Y}_{(n)} \left( \bigodot_{k \neq n} \mathbf{A}^{(k)} \right) \left( \bigcircledast_{k \neq n} \mathbf{A}^{(k)T} \mathbf{A}^{(k)} \right)^{\dagger}, \quad (n = 1, 2, \ldots, N), \tag{2.14}$$

where "$\dagger$" denotes the pseudo-inverse. A fast implementation of ALS for 3-way tensor [33] reduces the expensive computation of $\mathbf{Y}_{(n)} \left( \bigodot_{k \neq n} \mathbf{A}^{(k)} \right)$. Unfortunately, this algorithm cannot be generalized to higher orders [34]. The all-at-once algorithms such as OPT [1], PMF3, the damped Gauss-Newton (dGN) algorithms [24–26, 32, 33] compute gradients in their update rules

$$\mathfrak{a} \leftarrow \mathfrak{a} - \eta \, \mathbf{g}, \quad \eta > 0, \tag{2.15}$$

or

$$\mathfrak{a} \leftarrow \mathfrak{a} - (\mathbf{H} + \mu \mathbf{I}_{RT})^{-1} \, \mathbf{g}, \quad \mu > 0, T = \sum_n I_n, \tag{2.16}$$

where $\mathfrak{a} = \left[ \text{vec}\left(\mathbf{A}^{(1)}\right)^T \cdots \text{vec}\left(\mathbf{A}^{(n)}\right)^T \cdots \text{vec}\left(\mathbf{A}^{(N)}\right)^T \right]^T$, $\mathbf{H}$ denotes the (approximate) Hessian and $\mathbf{g}$ is the gradient defined as

$$\mathbf{g} = \frac{\partial D}{\partial \mathfrak{a}} = \left[ \left( \frac{\partial D}{\partial \text{vec}(\mathbf{A}^{(1)})} \right)^T \cdots \left( \frac{\partial D}{\partial \text{vec}(\mathbf{A}^{(n)})} \right)^T \cdots \left( \frac{\partial D}{\partial \text{vec}(\mathbf{A}^{(N)})} \right)^T \right]^T. \tag{2.17}$$

For a nonnegative tensor factorization, the well-known multiplicative algorithm [9, 23] also involves the CP gradients

$$\mathbf{A}^{(n)} \leftarrow \mathbf{A}^{(n)} \circledast \left( \mathbf{Y}_{(n)} \left( \bigodot_{k \neq n} \mathbf{A}^{(k)} \right) \right) \oslash \left( \mathbf{A}^{(n)} \left( \bigcircledast_{k \neq n} \mathbf{A}^{(k)T} \mathbf{A}^{(k)} \right) \right), \quad (n = 1, 2, \ldots, N). \tag{2.18}$$

The direct computation of the product $\mathbf{Y}_{(n)} \left( \bigodot_{k \neq n} \mathbf{A}^{(k)} \right)$ for single mode is illustrated in Algorithm 1, and is implemented in the `mttkrp` function of the Matlab Tensor toolbox [5, 6].

## 3. Fast Computation of CP Gradient.

### 3.1. Order of Dimensions.
The CP gradient $\mathbf{G}^{(n)}$ given by

$$\mathbf{G}^{(n)} = \mathbf{Y}_{(n)} \left(\underset{k \neq n}{\odot} \mathbf{A}^{(k)}\right) = \left[\mathbf{Y}_{(n)} \underset{k \neq n}{\otimes} \boldsymbol{a}_1^{(k)}, \; \mathbf{Y}_{(n)} \underset{k \neq n}{\otimes} \boldsymbol{a}_2^{(k)}, \ldots, \mathbf{Y}_{(n)} \underset{k \neq n}{\otimes} \boldsymbol{a}_R^{(k)}\right] \quad \in \mathbb{R}^{I_n \times R} \quad (3.1)$$

involves $R$ products

$$\boldsymbol{g}_r^{(n)} = \mathbf{Y}_{(n)} \left(\underset{k \neq n}{\otimes} \boldsymbol{a}_r^{(k)}\right) = \boldsymbol{\mathcal{Y}} \bar{\times}_1 \boldsymbol{a}_r^{(1)} \cdots \bar{\times}_{n-1} \boldsymbol{a}_r^{(n-1)} \bar{\times}_{n+1} \boldsymbol{a}_r^{(n+1)} \cdots \bar{\times}_N \boldsymbol{a}_r^{(N)}, \quad r = 1, \ldots R. \quad (3.2)$$

The Kronecker products $\boldsymbol{t} = \bigotimes_{k \neq n} \boldsymbol{a}_r^{(k)}$ can be efficiently computed by the following scheme [5,6]

$$\begin{aligned} \boldsymbol{t} \leftarrow \boldsymbol{a}^{(2)} \otimes \boldsymbol{a}^{(1)}, \quad & \boldsymbol{t} \leftarrow \boldsymbol{a}^{(3)} \otimes \boldsymbol{t}, \quad \ldots, \boldsymbol{t} \leftarrow \boldsymbol{a}^{(n-1)} \otimes \boldsymbol{t}, \\ \boldsymbol{t} \leftarrow \boldsymbol{a}^{(n+1)} \otimes \boldsymbol{t}, \quad & \boldsymbol{t} \leftarrow \boldsymbol{a}^{(n+2)} \otimes \boldsymbol{t}, \ldots, \boldsymbol{t} \leftarrow \boldsymbol{a}^{(N)} \otimes \boldsymbol{t}, \end{aligned} \quad (3.3)$$

with a computational cost of $O(\sum_{k=2}^{n-1} J_k + \frac{1}{I_n} \sum_{k=n+1}^{N} J_k)$.

Assuming that $I_N < I_{N-1}$, we transpose $\boldsymbol{\mathcal{Y}}$ following $\boldsymbol{p} = [1:N-2, N, N-1]$ to obtain tensor $\boldsymbol{\mathcal{Y}}^{<p>}$ of the size $I_1 \times \cdots \times I_{N-2} \times I_N \times I_{N-1}$. The tensor-vector products in (3.2) can be expressed by

$$\begin{aligned} \boldsymbol{\mathcal{Y}} \bar{\times}_{k=1}^{n-1} \{\boldsymbol{a}_r^{(k)}\} \bar{\times}_{l=n+1}^{N} \{\boldsymbol{a}_r^{(l)}\} &= \boldsymbol{\mathcal{Y}}^{<p>} \bar{\times}_{k=1}^{n-1} \{\boldsymbol{a}_r^{(k)}\} \bar{\times}_{l=n+1}^{N-2} \{\boldsymbol{a}_r^{(l)}\} \bar{\times}_{N-1} \boldsymbol{a}_r^{(N)} \bar{\times}_N \boldsymbol{a}_r^{(N-1)} \\ &= \mathbf{Y}_{(n)}^{<p>} \left(\boldsymbol{a}^{(N-1)} \otimes \boldsymbol{a}^{(N)} \otimes \left(\bigotimes_{l=n+1}^{N-2} \boldsymbol{a}_r^{(l)}\right) \otimes \left(\bigotimes_{k=1}^{n-1} \boldsymbol{a}_r^{(k)}\right)\right). \end{aligned} \quad (3.4)$$

According to the above computation scheme in (3.3), the Kronecker products in (3.4) require a computational cost of

$$O\left(\sum_{k=2}^{n-1} J_k + \frac{1}{I_n}\left(\sum_{k=n+1}^{N-2} J_k + J_{N-2} I_N + J_N\right)\right) < O\left(\sum_{k=2}^{n-1} J_k + \frac{1}{I_n} \sum_{k=n+1}^{N} J_k\right) \quad (3.5)$$

by noting that $J_{N-1} = J_{N-2} I_{N-1}$. As a result, in order to efficiently compute $\mathbf{G}^{(n)}$, we need to permute the tensor $\boldsymbol{\mathcal{Y}}$ such that $I_1 \leq I_2 \leq \cdots \leq I_N$. Hereinafter, we implicitly assume that the data tensor has been rearranged in the ascending order of its dimensions.

From (3.3), computation of $\mathbf{G}^{(n)}$ in (3.1) requires a number of multiplications of

$$M_{Alg.\,1}(n) = R\left(J_N + \sum_{k=2}^{n-1} J_k + \frac{1}{I_n} \sum_{k=n+1}^{N} J_k\right). \quad (3.6)$$

In a particular case when $I_n = I, \forall n$, Algorithm 1 executes a number of multiplications of $M_{Alg.\,1}(n) = R \sum_{k=2}^{N} I^k$.

### 3.2. Fast Gradient with Respect to A Specific Factor.
The direct computation of $\mathbf{G}^{(n)} = [\boldsymbol{g}_r^{(n)}]$ in (3.1) involves the tensor unfolding $\mathbf{Y}_{(n)}$ which is relatively slow to obtain for $1 < n < N$, due to accessing non-contiguous blocks of entries. We note that vectors $\boldsymbol{g}_r^{(n)}$ can be expressed in an equivalent form consisting of tensor-vector products $\boldsymbol{\mathcal{Y}} \bar{\times}_{k=1}^{n-1} \boldsymbol{a}_r^{(k)}$ and $\boldsymbol{\mathcal{Y}} \bar{\times}_{l=n+1}^{N} \boldsymbol{a}_r^{(l)}$ on the right side and left side of $n$, that is

$$\boldsymbol{g}_r^{(n)} = \mathbf{Y}_{(n)} \left(\bigotimes_{k \neq n} \boldsymbol{a}_r^{(k)}\right) = \left(\boldsymbol{\mathcal{Y}} \bar{\times}_{k=1}^{n-1} \boldsymbol{a}_r^{(k)}\right) \bar{\times}_{l=2}^{N-n+1} \boldsymbol{a}_r^{(l+n-1)}, \quad (3.7)$$

or

$$g_r^{(n)} = \left(\mathcal{Y} \bar{\times}_{l=n+1}^{N} a_r^{(l)}\right) \bar{\times}_{k=1}^{n-1} a_r^{(k)}. \tag{3.8}$$

We show in the sequel, that the former way, (3.7), is less computationally demanding for $J_n \leq K_{n-1}$, and the latter way, (3.8) is less demanding in the opposite case, $J_n > K_{n-1}$.

Note that the inner tensor-vector products in (3.7) and (3.8) can be efficiently computed through $\mathbf{Y}_{(1:n-1)}$ and $\mathbf{Y}_{(1:n)}$ as

$$\mathcal{L}^{(r,n)} \triangleq \mathcal{Y} \bar{\times}_{k=1}^{n-1} a_r^{(k)} = \text{reshape}\left(\mathbf{Y}_{(1:n-1)}^T \bigotimes_{k=1}^{n-1} a_r^{(k)}, [I_n, \ldots, I_N]\right), \tag{3.9}$$

$$\mathcal{R}^{(r,n)} \triangleq \mathcal{Y} \bar{\times}_{k=n+1}^{N} a_r^{(k)} = \text{reshape}\left(\mathbf{Y}_{(1:n)} \bigotimes_{k=n+1}^{N} a_r^{(k)}, [I_1, \ldots, I_n]\right). \tag{3.10}$$

It means that the reshapings of $\mathcal{Y}$ to $\mathbf{Y}_{(n)}$ with $1 < n < N$ are avoided.

Let us discuss complexity of the two ways of computing the gradients $g_r^{(n)}$ separately.

**3.2.1. The case $J_n \leq K_{n-1}$.** The computation proceeds first by computing the $n$-dimensional tensors $\mathcal{R}^{(r,n)}$ defined in (3.10) for all $r = 1, \ldots, R$. This operation requires the number of multiplications of

$$M_{RL1} = R\left(J_N + \frac{1}{J_n} \sum_{k=n+2}^{N} J_k\right). \tag{3.11}$$

The second step consists in computing $g_r^{(n)}$ as

$$g_r^{(n)} = \mathbf{R}_{(n)}^{(r,n)} \left(\bigotimes_{k=1}^{n-1} a_r^{(k)}\right). \tag{3.12}$$

where $\mathbf{R}_{(n)}^{(r,n)}$ is the mode-$n$ unfolding of $\mathcal{R}^{(r,n)}$. The second step has the complexity

$$M_{RL2} = R\left(J_n + \sum_{k=2}^{n-1} J_k\right). \tag{3.13}$$

The total number of multiplications is

$$M_{RL}(n) = R\left(J_N + \frac{1}{J_n} \sum_{k=n+2}^{N} J_k + \sum_{k=2}^{n} J_k\right) \leq R\left(J_N + \frac{1}{I_n} \sum_{k=n+2}^{N} J_k + J_{n-1}I_{n+1} + \sum_{k=2}^{n-1} J_k\right)$$

$$= R\left(J_N + \frac{1}{I_n} \sum_{k=n+1}^{N} J_k + \sum_{k=2}^{n-1} J_k\right) = M_{Alg.\ 1}(n), \quad 1 < n < N, \tag{3.14}$$

which is less than that of Algorithm 1 due to $J_n > I_n$ and $I_n \leq I_{n+1}$. That means the right-to-left projections should be faster than Algorithm 1.

**3.2.2. The case $J_n > K_{n-1}$.** Here we build up $(N - n + 1)$-D tensors $\mathcal{L}^{(r,n)}$ of size $I_n \times I_{n+1} \times \cdots \times I_N$ defined in (3.9). This step requires the number of multiplications of

$$M_{LR1} = R\left(J_N + \sum_{k=2}^{n-1} J_k\right). \tag{3.15}$$

The second step consists in computing the product

$$g_r = \mathbf{L}_{(1)}^{(r,n)} \left( \bigotimes_{k=n+1}^{N} a_r^{(k)} \right), \tag{3.16}$$

where $\mathbf{L}_{(1)}^{(r,n)}$ is mode-1 unfolding of $\mathcal{L}^{(r,n)}$. The number of multiplications in (3.16) is given by

$$M_{LR2} = R \left( K_{n-1} + \frac{1}{J_n} \sum_{k=n+2}^{N} J_k \right). \tag{3.17}$$

From (3.15) and (3.17), the proposed algorithm requires a total number of multiplications of

$$M_{LR}(n) = R \left( J_N + \sum_{k=2}^{n-1} J_k + K_{n-1} + \frac{1}{J_n} \sum_{k=n+2}^{N} J_k \right) < R \left( J_N + \sum_{k=2}^{n-1} J_k + J_n + \frac{1}{J_n} \sum_{k=n+2}^{N} J_k \right)$$
$$= M_{RL}(n) < M_{Alg.\ 1}(n), \quad 1 < n < N, \tag{3.18}$$

which is less than $M_{RL}$ in (3.14) in the previous subsection and of Algorithm 1 due to $K_{n-1} < J_n$ and $J_n > I_n$. That means the left-to-right projections should be faster than Algorithm 1.

**3.3. Fast CP gradient From Adjacent Ones.** CP algorithms available in the literature compute all $\mathbf{G}^{(n)}$ $n = 1, 2, \ldots, N$ either sequentially (in alternating algorithms [2, 8, 16, 18, 23, 30, 31]) or simultaneously (as in all-at-once algorithms [1, 24–26, 32, 33], line-search [27, 33]). This section will present a fast method to compute the gradients recursively for all $n = 1, \ldots, N$.

Note that

$$\mathcal{R}^{(r,n)} = \mathcal{Y} \bar{\times}_{k=n+1}^{N} a_r^{(k)} = \mathcal{R}^{(r,n+1)} \bar{\times}_{n+1} a_r^{(n+1)}, \tag{3.19}$$

or

$$\text{vec}\left(\mathcal{R}^{(r,n)}\right) = \mathbf{R}_{(1:n)}^{(r,n+1)} a_r^{(n+1)}. \tag{3.20}$$

Similarly,

$$\mathcal{L}^{(r,n)} = \mathcal{Y} \bar{\times}_{k=1}^{n-1} a_r^{(k)} = \mathcal{L}^{(r,n-1)} \bar{\times}_1 a_r^{(n-1)}, \tag{3.21}$$

or

$$\text{vec}\left(\mathcal{L}^{(r,n)}\right) = \mathbf{L}_{(1)}^{(r,n-1)T} a_r^{(n-1)}. \tag{3.22}$$

By exploiting relations in (3.20) and (3.22), we can quickly derive $\mathcal{R}^{(n)}$ from $\mathcal{R}^{(n+1)}$ or $\mathcal{L}^{(n)}$ from $\mathcal{L}^{(n-1)}$ instead of fully computing them as in (3.10) and (3.9), respectively. The total number of multiplications of the algorithm is summarized in Table 4.1 and is lower than that of Algorithm 1.

The proposed algorithm to compute CP gradients over all modes is summarized as Algorithm 2. Gradient $\mathbf{G}^{(n^\star)}$ or $\mathbf{G}^{(n^{\star+1})}$ is first computed where $n^\star = \max\{n; J_n \leq K_n\}$. Gradients $\mathbf{G}^{(n)}$ for $n = n^\star - 1, n^\star - 2, \ldots, 1$, and $\mathbf{G}^{(n)}$ for $n = n^\star + 1, n^\star + 2, \ldots, N$ are then sequentially computed. Note that in addition to the lower number of multiplications, Algorithm 2 avoids unfolding $\mathbf{Y}_{(n)}$ ($1 < n < N$) which is time consuming. Therefore, the higher the dimension of the tensor is, the more significant the computational saving of Algorithm 2 with respect to Algorithm 1 is.

**Algorithm 2:** Fast Computation of $\mathbf{Y}_{(n)} \left( \bigodot_{k \neq n} \mathbf{A}^{(k)} \right)$ over all modes

**Input**: $\mathcal{Y}$: $(I_1 \times I_2 \times \cdots \times I_N)$, $N$ factors $\mathbf{A}^{(n)} \in \mathbb{R}^{I_n \times R}$
**Output**: $\mathbf{G}^{(n)} = \mathbf{Y}_{(n)} \left( \bigodot_{k \neq n} \mathbf{A}^{(k)} \right)$: $I_n \times R$, $n = 1, 2, \ldots, N$

**begin**
1.    $n^\star = \max\{n; J_n \leq K_n\}$ where $J_n = I_1 I_2 \ldots I_n$, $K_n = I_{n+1} \cdots I_N$
2.    **for** $n = n^\star, n^\star - 1, \ldots, 1, n^\star + 1, n^\star + 2, \ldots, N$ **do**
3.      **if** $(n == n^\star)$ **then**
4.        $\mathbf{R}^{(n)}_{(n)} = \texttt{reshape}(\mathcal{Y}, [J_n, K_n]) \left( \bigodot_{k=n+1}^{N} \mathbf{A}^{(k)} \right)$     % mode-$n$ unfolding of $\mathcal{R}^{(n)}$
5.        $\mathbf{G}^{(n)} = \texttt{cp\_gradient}(\mathcal{R}^{(n)}, \{\mathbf{A}\}, \textit{'right'})$
6.      **else if** $(n == n^\star + 1)$ **then**
7.        $\mathbf{L}^{(n)}_{(1)} = \left( \bigodot_{k=1}^{n-1} \mathbf{A}^{(k)} \right)^T \texttt{reshape}(\mathcal{Y}, [J_{n-1}, K_{n-1}])$     % mode-1 unfolding of $\mathcal{L}^{(n)}$
8.        $\mathbf{G}^{(n)} = \texttt{cp\_gradient}(\mathcal{L}^{(n)}, \{\mathbf{A}\}, \textit{'left'})$
9.      **else if** $(n < n^\star)$ **then**
10.        **for** $r = 1, 2, \ldots, R$ **do**       % Compute $\mathcal{R}^{(r,n)}$ as in (3.20)
11.          $\texttt{vec}(\mathcal{R}^{(r,n)}) \leftarrow \texttt{reshape}(\mathcal{R}^{(r,n+1)}, [J_n, I_{n+1}]) \, \boldsymbol{a}^{(n+1)}_r$
12.        $\mathbf{G}^{(n)} = \texttt{cp\_gradient}(\mathcal{R}^{(n)}, \{\mathbf{A}\}, \textit{'right'})$
13.      **else**
14.        **for** $r = 1, 2, \ldots, R$ **do**       % Compute $\mathcal{L}^{(r,n)}$ as in (3.22)
15.          $\texttt{vec}(\mathcal{L}^{(r,n)}) \leftarrow \texttt{reshape}(\mathcal{L}^{(r,n-1)}, [I_{n-1}, K_{n-1}])^T \, \boldsymbol{a}^{(n-1)}_r$
16.        $\mathbf{G}^{(n)} = \texttt{cp\_gradient}(\mathcal{L}^{(n)}, \{\mathbf{A}\}, \textit{'left'})$

**function** $\mathbf{G}^{(n)} = \texttt{cp\_gradient}(\mathcal{Z}^{(n)}, \{\mathbf{A}\}, \textit{side})$
17. **for** $r = 1, 2, \ldots, R$ **do**
18.    **switch** *side*
19.      **case** *'right'*: $\boldsymbol{g}^{(n)}_r = \texttt{reshape}(\mathcal{Z}^{(r,n)}, [J_{n-1}, I_n])^T \left( \bigotimes_{k=1}^{n-1} \boldsymbol{a}^{(k)}_r \right)$    % $\mathcal{Z}^{(r,n)} \equiv \mathcal{R}^{(r,n)}$ % in (3.12)
20.      **case** *'left'*: $\boldsymbol{g}^{(n)}_r = \texttt{reshape}(\mathcal{Z}^{(r,n)}, [I_n, K_n]) \left( \bigotimes_{k=n+1}^{N} \boldsymbol{a}^{(k)}_r \right)$    % $\mathcal{Z}^{(r,n)} \equiv \mathcal{L}^{(r,n)}$ % in (3.16)

$\mathcal{R}^{(n)}(i_1, \ldots, i_n, r) = \mathcal{R}^{(r,n)}(i_1, \ldots, i_n)$ in Step 4
$\mathcal{L}^{(n)}(r, i_n, \ldots, i_N) = \mathcal{L}^{(r,n)}(i_n, \ldots, i_N)$ in Step 7

**4. Fast ALS Algorithm.** This section presents a fast implementation of the CP ALS algorithm (2.14) in which gradients are computed using Algorithm 2. That is, the fast ALS algorithm is proposed to first update $\mathbf{A}^{(n^\star)}$ or $\mathbf{A}^{(n^\star+1)}$ instead of $\mathbf{A}^{(1)}$. The algorithm then updates sequentially $\mathbf{A}^{(n)}$ for $n = n^\star - 1, n^\star - 2, \ldots, 1$, and $\mathbf{A}^{(n)}$ for $n = n^\star + 1, n^\star + 2, \ldots, N$. The alternating update rules (2.14) are inserted in the "for" loop in Algorithm 2, and are executed after computing gradients $\mathbf{G}^{(n)}$. Such strategy requires a computational cost of order $O(RJ_N + NR^3)$ to complete updating all $\mathbf{A}^{(n)}$. Other alternative algorithms [9, 18, 23] can be accelerated in a similar way.

TABLE 4.1
*Comparison of the number of multiplications executed in methods to compute $\mathbf{Y}_{(n)} \left( \bigodot_{k \neq n} \mathbf{A}^{(k)} \right)$.*

|  | Number of multiplications | Unfoldings | - Order of Entries |
|---|---|---|---|
| $M_{Alg.\ 1}(n) = R\left(\sum_{k=2}^{n-1} J_k + \frac{1}{I_n}\sum_{k=n+1}^{N} J_k + J_N\right)$ | | $\mathbf{Y}_{(n)}$ | - change |
| $M_{Alg.\ 2}(n) = \begin{cases} R\left(\sum_{k=2}^{n-1} J_k + \min(J_n, K_{n-1}) + \frac{1}{J_n}\sum_{k=n+2}^{N} J_k + J_N\right), \\ \qquad\qquad n = n^\star, n^\star + 1, \\ R\sum_{k=2}^{n+1} J_k, \quad n < n^\star, \\ R\left(K_{n-2} + K_{n-1} + \sum_{k=n+2}^{N} J_k\right), \quad n > n^\star + 1 \end{cases}$ | | $\mathbf{Y}_{(1:n^\star)}$ or $\mathbf{Y}_{(1:n^\star-1)}$, $\mathbf{R}^{(r,n+1)}_{(1:n)}$ and $\mathbf{R}^{(r,n)}_{(n)}$ or $\mathbf{L}^{(r,n)}_{(1)}$ | - no-change |

**5. Simulations.** In order to verify the fast CP gradients (Algorithm 2), we compared the fast CP_ALS algorithm in Section 4 with the ordinary CP_ALS algorithm [8, 16] which was implemented in the Matlab Tensor toolbox [6] and used the direct computation of CP gradients in Algorithm 1 (`mttkrp` of the Tensor toolbox (ver. 2.4) [5, 6]). The Matlab codes of the fast CP gradient and fast CP_ALS are available at the following link: http://www.bsp.brain.riken.jp/~phan/fastCPgradient.rar. Random tensors with various dimensions $N = 3, 4, 5, 6$ were randomly generated with different sizes $I_n = I = 10, 20, \ldots, \forall n$. Both algorithms factorized the same data tensors into various rank $R = 1, 10, 20, \ldots, I$ using the same initialization values and in 20 iterations. There was not any stopping criterion for both algorithms. Execution time for each algorithm was measured using the stopwatch command: "tic" "toc" of MATLAB release 2011a on a computing server which has 1.8 GHz i7 processor and 4 GB memory. The Tucker compression was not used in the simulations.

Speed ratio is defined as the ratio between execution times per iterations of CP_ALS and the fast CP_ALS

$$\rho = \frac{Execution\_time_{ALS}}{Execution\_time_{fastALS}}. \tag{5.1}$$

The final results were averaged over at least 200 iterations. Fig. 5.1 illustrates speed-up ratio per iteration (times) for factorization of 3-D and 4-D tensors with different sizes $I$ and ranks $R$, whereas the speed-up ratios per iterations for 5-D and 6-D tensors are given in Table 5.1. In an extra example for decomposition of 7-D tensors with $I_n = 10, \forall n$, and $R = 10$, the ordinary ALS algorithm took an average 1.844 seconds per iteration, while the fast ALS took only 0.044 second per iteration, and achieved an average speed-up ratio of 42.2 times. For some data which consists of collinear factors, such as bottleneck or swamps [10], the ALS algorithm could execute thousands of iterations. Hence, the ALS algorithm [8, 16] needs at least 1 hour to run over 2000 iterations. Whereas the fast ALS algorithm executes the same number of iterations only in 1 or 2 minutes. This indicates a huge potential benefit of our fast algorithm.

**6. Conclusions.** CP gradients are always the largest workload of order $O\left(NR\prod_{n=1}^{N} I_n\right)$ in most CP algorithms such as the ALS and all-at-one algorithms. Moreover, the computation can be time consuming due to unfoldings $\mathbf{Y}_{(n)}$ with $1 < n < N$. The fast computation method has been proposed to avoid $\mathbf{Y}_{(n)}$ for $1 < n < N$, and has an approximate computational cost of order $O\left(R \prod_{n=1}^{N} I_n\right)$. Especially, the fast CP gradients can be used to accelerate any CP and NTF algorithms. Experimental results show that our algorithm can be about 8 times faster than the direct computation for 4-D tensors, and it can be up to 20-30 times for higher dimensional tensors.

TABLE 5.1
*Comparison execution times (seconds) between the CP ALS algorithm and the fast CP ALS using the fast CP gradients for random data tensors of size $I_1 = I_2 = \ldots = I_N = I$. Results for each combination (N, I, R) consist of execution times and speed-up ratio between two algorithms as indicated in the below minitab. The results were averaged over at least 200 runs on a computer which has 1.8 GHz i7 processor and 4 GB memory.*

| N, I | R |  |  |  |  |  |  |  |  |
|---|---|---|---|---|---|---|---|---|---|
|  | 1 |  | 10 |  | 20 |  | 30 |  | 40 |
| 5, 20 | $\frac{9\,10^{-3}}{3\,10^{-4}}$ | 33.1 | $\frac{0.02}{9\,10^{-4}}$ | 19.3 | $\frac{0.02}{1.2\,10^{-3}}$ | 19.6 |  |  |  |  |
| 5, 30 | $\frac{0.05}{9\,10^{-4}}$ | 62.1 | $\frac{0.09}{4.3\,10^{-3}}$ | 21.9 | $\frac{0.1}{6.2\,10^{-3}}$ | 21.2 | $\frac{0.2}{0.011}$ | 17.3 |  |  |
| 5, 40 | $\frac{0.3}{3\,10^{-3}}$ | 83.1 | $\frac{0.35}{0.016}$ | 21.6 | $\frac{0.46}{0.023}$ | 20.1 | $\frac{0.63}{0.036}$ | 17.6 | $\frac{1.4}{0.046}$ | 31.6 |
| 6, 20 | $\frac{0.2}{2\,10^{-3}}$ | 102.9 | $\frac{0.34}{0.01}$ | 24.2 | $\frac{0.5}{0.015}$ | 34.7 |  |  |  |  |

| $\frac{Execution\_time_{ALS} \text{ (seconds)}}{Execution\_time_{fastALS} \text{ (seconds)}}$ | ratio |
|---|---|

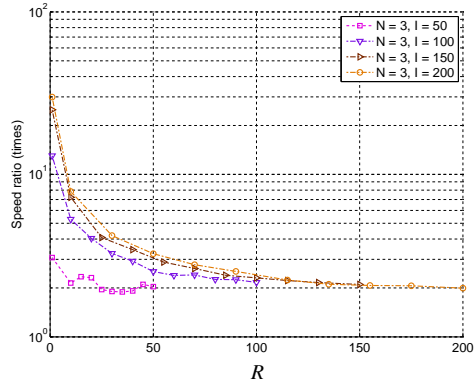 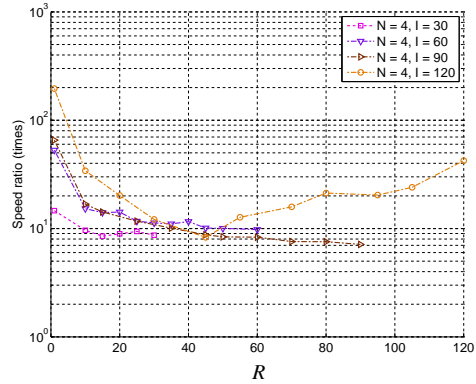

(a) 3-D tensors  (b) 4-D tensors

FIG. 5.1. *Speed-up ratios per iteration (in logarithmic scale) between the ordinary CP_ALS algorithm and its fast implementation using the fast CP gradients for factorization of 3-D and 4-D tensors with various sizes $I_n = I, \forall n$, and ranks R.*